\documentclass[preprint,12pt,sort&compress]{elsarticle}

\usepackage{amssymb}
\usepackage{amsthm}
\usepackage{amsmath}
\usepackage{mathrsfs}
\usepackage{graphicx}
\usepackage{epstopdf}
\usepackage{float}
\usepackage{caption}
\usepackage{subcaption}
\usepackage{bm}
\usepackage{bbm}
\usepackage{mathrsfs}
\usepackage{cleveref}
\usepackage{soul}
\usepackage{accents}
\usepackage{graphicx}
\usepackage{courier} % To write in Courier font (to write code basically)
\usepackage{color} % To write in different colors
\journal{Computational Mechanics}

\makeatletter
\def\@author#1{\g@addto@macro\elsauthors{\normalsize%
    \def\baselinestretch{1}%
    \upshape\authorsep#1\unskip\textsuperscript{%
      \ifx\@fnmark\@empty\else\unskip\sep\@fnmark\let\sep=,\fi
      \ifx\@corref\@empty\else\unskip\sep\@corref\let\sep=,\fi
      }%
    \def\authorsep{\unskip,\space}%
    \global\let\@fnmark\@empty
    \global\let\@corref\@empty  %% Added
    \global\let\sep\@empty}%
    \@eadauthor={#1}
}
\makeatother

\begin{document}

\begin{frontmatter}



\title{Gradient plasticity crack tip characterization by means of the extended finite element method}


\author{E. Mart\'{\i}nez-Pa\~neda\corref{cor1}\fnref{Uniovi}}
\ead{mail@empaneda.com}

\author{S. Natarajan\fnref{IIT}}

\author{S. Bordas\fnref{UL,UC}}

\address[Uniovi]{Department of Mechanical Engineering, Solid Mechanics, Technical University of Denmark, DK-2800 Kgs. Lyngby, Denmark}

\address[IIT]{Department of Mechanical Engineering, Indian Institute of Technology, Madras, Chennai 600036, India}

\address[UL]{University of Luxembourg, Research Unit of Engineering Science, Campus Kirchberg, Luxembourg}

\address[UC]{Cardiff University, School of Engineering, Cardiff CF24 3AA, Wales, UK}

\cortext[cor1]{Corresponding author.}

\begin{abstract}
Strain gradient plasticity theories are being widely used for fracture assessment, as they provide a richer description of crack tip fields by incorporating the influence of geometrically necessary dislocations. Characterizing the behavior at the small scales involved in crack tip deformation requires, however, the use of a very refined mesh within microns to the crack. In this work a novel and efficient gradient-enhanced numerical framework is developed by means of the extended finite element method (X-FEM). A mechanism-based gradient plasticity model is employed and the approximation of the displacement field is enriched with the stress singularity of the gradient-dominated solution. Results reveal that the proposed numerical methodology largely outperforms the standard finite element approach. The present work could have important implications on the use of microstructurally-motivated models in large scale applications. The non-linear X-FEM code developed in MATLAB can be downloaded from www.empaneda.com/codes
\end{abstract}

\begin{keyword}

Strain gradient plasticity \sep Extended Finite Element Method \sep Crack tip fields \sep material length scale \sep MATLAB



\end{keyword}

\end{frontmatter}


\section{Introduction}
\label{Introduction}

Experiments have consistently shown that metallic materials display strong size effects at the micron scale, with smaller being harder. As a result, a significant body of research has been devoted to model this size dependent plastic phenomenon (see, e.g., \cite{DF13,BB16,M16} and references therein). At the continuum level, phenomenological strain gradient plasticity (SGP) formulations have been developed to extend plasticity theory to small scales. Grounded on the physical notion of geometrically necessary dislocations (GNDs, associated with non-uniform plastic deformation), SGP theories relate the plastic work to both strains and strain gradients, introducing a length scale in the constitutive equations. Isotropic SGP formulations can be classified according to different criteria, one distinguishing between phenomenological theories \cite{A84,FH01} and microstructurally or mechanism-based ones \citep{G99,H04}. All these models aim at predicting the strengthening effects associated with dislocation interactions in an average sense, as opposed to the more refined characterization of explicit multiscale approaches \cite{H04b,V07,M08,M08b,O15}.\\

While growing interest in micro-technology motivated the development of SGP models at first, the influence of GNDs extends beyond micron-scale applications, as strains vary over microns in a wide range of engineering designs. Particularly, gradient-enhanced modeling of fracture and damage appears imperative - independently of the size of the specimen - as the plastic zone adjacent to the crack tip is physically small and contains strong spatial gradients of deformation. The experimental observation of cleavage fracture in the presence of significant plastic flow \cite{K02} has fostered significant interest in the role of the plastic strain gradient in fracture and damage assessment. Studies conducted in the framework of phenomenological \cite{K08,TN08,MG09} and mechanism-based theories \cite{J01,PY11,PY11b} have shown that GNDs near the crack tip promote local strain hardening and lead to a much higher stress level as compared with classic plasticity predictions. Very recently, Mart\'{\i}nez-Pa\~neda and co-workers \cite{MB15,MN16} have identified and quantified the relation between material parameters and the physical length over which gradient effects prominently enhance crack tip stresses. Their results have revealed the important influence of strain gradients on a wide range of fracture problems, being particularly relevant in hydrogen assisted cracking modeling due to the central role that the stress field close to the crack tip plays on both hydrogen diffusion and interface decohesion \cite{M16a,M16b}.\\

However, a comprehensive embrace of SGP theories has been hindered by the complexities associated with their numerical implementation. An appropriate characterization of gradient effects ahead of a crack requires the use of extremely refined meshes, with a characteristic element length of a few nanometers in the vicinity of the crack. The vast majority of fracture and damage studies under SGP theories (see, e.g., \cite{K08,TN08,MG09,J01,MB15,MN16}) have been conducted in the framework of the finite element method (FEM). \citet{PY11,PY11b} used the element-free Galerkin method with the aim of avoiding the lack of convergence associated with finite element schemes (especially as the element distortion becomes large or elements lose bearing capacity) at the expense of increasing the computational cost. In this work a novel numerical framework is proposed for crack tip assessment within strain gradient plasticity. The mechanism-based strain gradient (MSG) plasticity theory is adopted as a material model, being this choice motivated by the work by Shi \textit{et al.} \cite{S01}, who characterized - by means of combined analytical and numerical (Runge-Kutta) procedure - the stress-dominated asymptotic field around a crack tip within the aforementioned constitutive framework. The extended finite element method (X-FEM) can be therefore employed to enrich the solution significantly alleviating the degree of mesh refinement. In the present work, a novel non-linear X-FEM scheme is presented, which includes (i) gradient-enhanced asymptotic functions, (ii) linear and quadratic elements, (iii) a linear weighting function for the blending elements, (iv) an iterative solver for nonlinear systems and (v) an appropriate triangular integration scheme. Several numerical examples are addressed to illustrate the performance of the present numerical approach.

\section{Mechanism-based strain gradient (MSG) plasticity}
\label{Material model}

The theory of mechanism-based strain gradient plasticity \citep{G99,H04} is based on the Taylor dislocation model and therefore the shear flow stress $\tau$ is formulated in terms of the dislocation density $\rho$ as
\begin{equation}\label{Eq1MSG}
\tau = \alpha \mu b \sqrt{\rho}
\end{equation}

Here, $\mu$ is the shear modulus, $b$ is the magnitude of the Burgers vector and $\alpha$ is an empirical coefficient which takes values between 0.3 and 0.5. The dislocation density is composed of the sum of the density $\rho_S$ for statistically stored dislocations (SSDs) and the density $\rho_G$ for geometrically necessary dislocations (GNDs) as
\begin{equation}\label{Eq2MSG}
\rho = \rho_S + \rho_G
\end{equation}

The GND density $\rho_G$ is related to the effective plastic strain gradient $\eta^{p}$ by: 
\begin{equation}\label{Eq3MSG}
\rho_G = \overline{r}\frac{\eta^{p}}{b}
\end{equation}

\noindent where $\overline{r}$ is the Nye-factor which is assumed to be 1.90 for face-centered-cubic (fcc) polycrystals. The tensile flow stress $\sigma_{flow}$ is related to the shear flow stress $\tau$ by:
\begin{equation}\label{Eq4MSG}
\sigma_{flow} =M\tau
\end{equation}

\noindent with $M$ being the Taylor factor, taken to be 3.06 for fcc metals. Rearranging Eqs. (\ref{Eq1MSG}-\ref{Eq4MSG}) yields
\begin{equation}\label{Eq5MSG}
\sigma_{flow} =M\alpha \mu b \sqrt{\rho_{S}+\overline{r}\frac{\eta^{p}}{b}}
\end{equation}

The SSD density $\rho_{S}$ can be determined from (\ref{Eq5MSG}) knowing the relation in uniaxial tension between the flow stress and the material stress-strain curve as follows
\begin{equation}\label{Eq6MSG}
\rho_{S} = [\sigma_{ref}f(\varepsilon^{p})/(M\alpha \mu b)]^2
\end{equation}

Here $\sigma_{ref}$ is a reference stress and $f$ is a nondimensional function of the plastic strain $\varepsilon^{p}$ determined from the uniaxial stress-strain curve. Substituting back into (\ref{Eq5MSG}), $\sigma_{flow}$ yields:
\begin{equation}\label{EqSflow}
\sigma_{flow} =\sigma_{ref} \sqrt{f^2(\varepsilon^{p})+l\eta^{p}}
\end{equation}

\noindent where $l$ is the intrinsic material length based on parameters of elasticity ($\mu$), plasticity ($\sigma_{ref}$) and atomic spacing ($b$):
\begin{equation}\label{Eqell}
l=M^2\overline{r}\alpha^2 \left(\frac{\mu}{\sigma_{ref}}\right)^2b=18\alpha^2\left(\frac{\mu}{\sigma_{ref}}\right)^2b
\end{equation}

\citet{G99} used three quadratic invariants of the plastic strain gradient tensor to represent the effective plastic strain gradient $\eta^{p}$ as
\begin{equation}
\eta^{p}=\sqrt{c_1 \eta^{p}_{iik} \eta^{p}_{jjk} + c_2 \eta^{p}_{ijk} \eta^{p}_{ijk} + c_3 \eta^{p}_{ijk} \eta^{p}_{kji}}
\end{equation}

The coefficients were determined to be equal to $c_1=0$, $c_2=1/4$ and $c_3=0$ from three dislocation models for bending, torsion and void growth, leading to
\begin{equation}\label{Eq:grad}
\eta^{p}=\sqrt{\frac{1}{4}\eta^{p}_{ijk} \eta^{p}_{ijk}}
\end{equation}

\noindent where the components of the strain gradient tensor are obtained by $\eta^{p}_{ijk}=\varepsilon^{p}_{ik,j}+\varepsilon^{p}_{jk,i}-\varepsilon^{p}_{ij,k}$.\\

As it is based on the Taylor dislocation model, which represents an average of dislocation activities, the MSG plasticity theory is only applicable at a scale much larger than the average dislocation spacing. For common values of dislocation density in metals, the lower limit of physical validity of the SGP theories based on Taylor’s dislocation model is approximately 100 nm.\\

Shi and co-workers \cite{S01} characterized the stress-dominated asymptotic field around a mode I crack tip in MSG plasticity by solving iteratively through Runge-Kutta a fifth order Ordinary Differential Equation (ODE). The numerical shooting method was employed to enforce two crack-face stress-traction free conditions and subsequently obtain the power of the stress singularity, roughly $r^{-2/3}$. A similar result was obtained through finite element (FE) analysis by Jiang \textit{et al.} \cite{J01}. The power of the stress singularity in MSG plasticity is therefore independent of the strain hardening exponent $n$. This is due to the fact that the strain gradient becomes more singular than the strain near the crack tip and dominates the contribution to the flow stress in (\ref{EqSflow}). From a physical viewpoint, this indicates that the density of GNDs $\rho_G$ in the vicinity of the crack tip is significantly larger than the density of SSDs $\rho_S$.\\

\begin{figure}[H]
\centering
\includegraphics[scale=0.5]{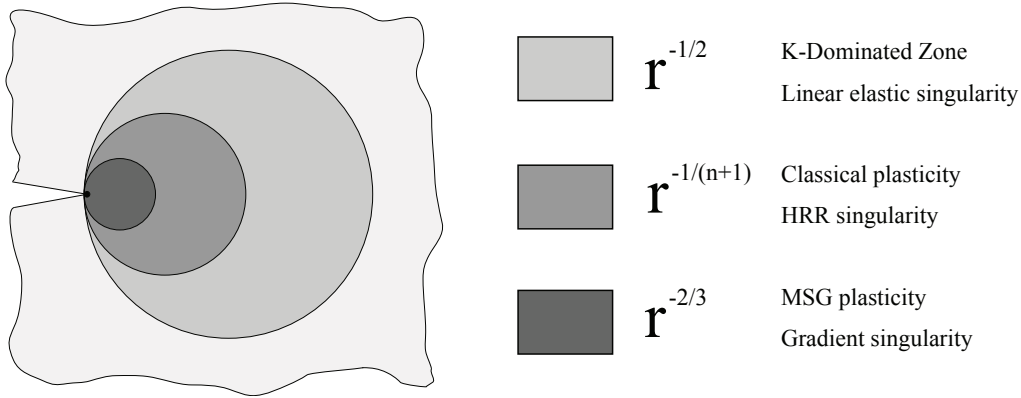}
\caption{Schematic diagram of the different domains surrounding the crack tip. Three regions are identified as a function of asymptotic stress fields: the linear elastic solution, the HRR solution and the MSG plasticity solution.}
\label{fig:Domain}
\end{figure}

Consequently, crack tip fields can be divided in several domains, as depicted in Fig. \ref{fig:Domain}. Far away from the crack tip deformation is elastic and the asymptotic stress field is governed by the linear elastic singularity. When the effective stress overcomes the initial yield stress $\sigma_Y$, plastic deformations occur and the stress field is characterized by the Hutchinson, Rice and Rosengren (HRR) \cite{H68,RR68} solution. As the distance to the crack tip decreases to the order of a few microns, large gradients of plastic strain promote dislocation hardening and the stress field is described by the asymptotic stress singularity of MSG plasticity.\\

\section{Numerical framework}
\label{X-FEM and FEM formulation}

\subsection{Finite element method}

As a function of their order, two different classes of SGP theories can be identified. One involves higher order stresses and therefore requires extra boundary conditions; the other does not involve higher order terms, and gradient effects come into play via the incremental plastic moduli. With the aim of employing mechanism-based SGP formulations within a lower order setup, Huang \textit{et al.} \cite{H04} developed what is referred to as the Conventional Mechanism-based Strain Gradient (CMSG) plasticity theory. It is also based on Taylor's dislocation model (i.e., MSG plasticity), but it does not involve higher order terms and therefore falls into the SGP framework that preserves the structure of classic plasticity. Consequently, the plastic strain gradient appears only in the constitutive model, and the equilibrium equations
and boundary conditions are the same as in conventional continuum
theories. This lower order scheme is adopted in the present work to characterize gradient effects from a mechanism-based approach, as it does not suffer convergence problems when addressing numerically demanding problems, such as crack tip deformation, unlike its higher order counterpart (see \cite{H03,MB15}). While higher order formulations are necessarily needed to model constraints on dislocation movement (see \cite{NH03b}), in MSG plasticity the differences between the higher order and the lower order
versions are restricted to a very thin boundary layer ($\approx$ 10 nm) \cite{S01,Q04}. Consequently, higher order boundary conditions essentially have no effect on the stress distribution at a distance of more than 10 nm away from the crack tip - well below its lower limit of physical validity -  and CMSG plasticity predicts exactly the same results as its higher order counterpart in the region of interest.\\

To avoid the use of higher order stresses, \citet{H04} used a viscoplastic formulation where the plastic strain rate $\dot{\varepsilon}^{p}$ is given in terms of the effective stress $\sigma_e$ rather than its rate $\dot{\sigma}_e$. The strain rate and time dependence is suppressed by adopting a viscoplastic power-law of the form,

\begin{equation}
\dot{\varepsilon}^{p} = \dot{\varepsilon} \left [\frac{\sigma_e}{\sigma_{ref} \sqrt{f^{2}(\varepsilon^{p})+l\eta^{p}}} \right]^{m}
\end{equation}

\noindent where the visco-plastic exponent is taken to fairly large values ($m\geq20$) in order to suppress rate effects. Taking into account that the volumetric and deviatoric strain rates are related to the stress rate in the same way as in conventional plasticity, the constitutive equation yields:

\begin{equation}
\dot{\sigma}_{ij}=K\dot{\varepsilon}_{kk}\delta_{ij}+2\mu \left\{\dot{\varepsilon}'_{ij} - \frac{3\dot{\varepsilon}}{2\sigma_e}\left[\frac{\sigma_e}{\sigma_{flow}} \right]^{m}\dot{\sigma}'_{ij} \right\}
\end{equation}

Since higher-order terms are not involved, the governing equations of CMSG plasticity are essentially the same as those in conventional plasticity and the FE implementation is relatively straightforward. The plastic strain gradient is obtained by numerical differentiation within the element through the shape functions. In order to do so, a surface is first created by linearly interpolating the incremental values of the plastic strains $\Delta \varepsilon^p_{ij}$ at the Gauss integration points in the entire model. Subsequently, the values of $\Delta \varepsilon^p_{ij}$ are sampled at the nodal locations. An almost identical procedure could be employed to map history-dependent variables in crack propagation studies and other cases where element subdivision takes place.\\

In order to validate the finite element implementation, crack tip fields are evaluated by means of a boundary layer formulation, where the crack region is contained by a circular zone and a Mode I load is applied at the remote circular boundary through a prescribed displacement: 

\begin{equation}
u(r,\theta)=K_I \frac{1+\nu}{E} \sqrt{\frac{r}{2\pi}}cos\left(\frac{\theta}{2}\right)(3-4\nu-cos\theta)
\end{equation}

\begin{equation}
v(r,\theta)=K_I \frac{1+\nu}{E} \sqrt{\frac{r}{2\pi}}sin\left(\frac{\theta}{2}\right)(3-4\nu-cos\theta)
\end{equation}

Here, $u$ and $v$ are the horizontal and vertical components of the displacement boundary condition, $r$ and $\theta$ the radial and angular coordinates in a polar coordinate system centered at the crack tip, $E$ is Young's modulus, $\nu$ is the Poisson ratio of the material and $K_I$ is the applied stress intensity factor, which quantifies the remote load. Plane strain conditions are assumed and only the upper half of the circular domain is modeled due to symmetry. A sufficiently large outer radius $R$ is defined and the entire specimen is discretized by means of eight-noded quadrilateral elements with reduced integration. Different mesh densities were used to study convergence behavior, with the typical number of elements being around 4000. With the aim of accurately characterizing the influence of the strain gradient a very refined mesh is used near the crack tip, where the size of the elements is on the order of very few nanometers. The following set of non-dimensional material parameters is considered:
\begin{equation}\label{Eq:prop}
n=5, \, \, \, \, \, \, \, \, \, \, \, \, \frac{\sigma_Y}{E}=0.2\%, \, \, \, \, \, \, \, \, \, \, \, \, \nu=0.3
\end{equation}

\noindent An isotropic power law material is adopted according to
\begin{equation}\label{Eq:PowerLaw}
\sigma=\sigma_Y \left( 1 + \frac{E \varepsilon^p}{\sigma_Y} \right)^{\left( \frac{1}{n} \right)}
\end{equation}

\begin{figure}[H]
\centering
\includegraphics[scale=0.75]{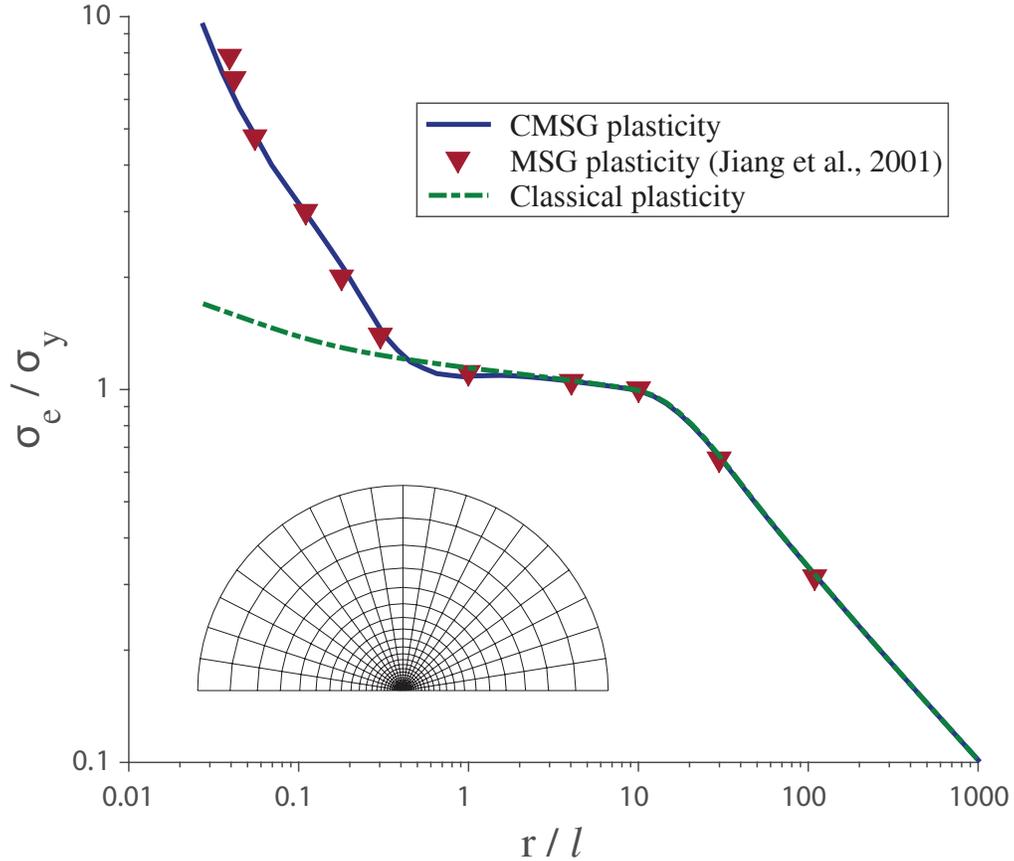}
\caption{Effective stress distribution ahead of the crack tip. Comparison between MSG plasticity predictions (symbols, taken from \cite{J01}), CMSG plasticity (solid line) and conventional $J_2$ plasticity (dashed line). The figure is in a double logarithm scale and $\sigma_e$ and $r$ are normalized by $\sigma_y$ and $l$ respectively. An external applied load of $K_{I}=20\sigma_Y \sqrt{l}$ is assumed and the following material properties are adopted: $\sigma_Y=0.2\%E$, $\nu=0.3$, $N=0.2$ and $l=3.53$ $\mu$m.}
\label{fig:Validation}
\end{figure}

The reference stress of (\ref{Eq6MSG}) will correspond to $\sigma_{ref}=\sigma_Y \left( \frac{E}{\sigma_Y}\right)^{\left( \frac{1}{n} \right)}$ and $f(\varepsilon^{p})=\left(\varepsilon^{p}+\frac{\sigma_Y}{E} \right)^{\left( \frac{1}{n} \right)}$. Fig. \ref{fig:Validation} shows, in a double logarithm diagram, the normalized effective stress $\sigma_e / \sigma_Y$ versus the normalized distance $r/l$ ahead of the crack tip ($\theta=1.014^{\circ}$) for an external applied load of $K_I=20 \sigma_Y \sqrt{l}$. As it can be seen in the figure, a very good agreement is obtained between the stress distributions obtained by means of the CMSG theory and MSG plasticity (taken from \cite{J01}); proving the suitability of CMSG plasticity in the present study, since higher order boundary conditions do not influence crack tip fields within its physical domain of validity.\\

As it can be seen in the figure, SGP predictions agree with classic plasticity away from the crack tip but become much larger within tens of microns from it. In agreement with the work by Shi \textit{et al.} \cite{S01}, the stress field in MSG plasticity is more singular than both the HRR field and the linear elastic K field. This GND-enhanced singularity can be incorporated within an X-FEM framework to avoid the use of extremely refined meshes and the numerical problems associated. 

\subsection{Extended Finite Element Method}
\label{X-FEM}

The approximation power of the FEM can be further enhanced by augmenting \emph{suitable} functions to the finite element space; these functions represent the local nature of the solution. This can be achieved by means of the X-FEM, a numerical enrichment strategy within the framework of the Partition of Unity Method (PUM). The displacement approximation can be thus decomposed into a standard part and an enriched part,
\begin{equation}
u_i^h =\underbrace{\sum_{I \in \mathcal{N}^{\rm fem}} N_i^I u_i^I}_\text{Standard} + \underbrace{\sum_{J \in \mathcal{N}^{\rm c}} N_i^J  H(\phi) a_i^J +  \sum_{K \in \mathcal{N}^{\rm f}} N_i^K \sum_{\alpha=1}^n F_{\alpha} (r, \theta ) b_i^{K \alpha}}_\text{Enriched}
\end{equation}

\noindent where $\mathcal{N}^{\rm fem}$ is the set of all nodes in the FE mesh, $\mathcal{N}^{\rm c}$ is the set of nodes whose shape function support is cut by the crack interior and $\mathcal{N}^{\rm f}$ is the set of nodes whose shape function support is cut by the crack tip. $H(\phi)$ and $F_{\alpha} (r, \theta )$ are the enrichment functions chosen to respectively capture the displacement jump across the crack surface and the singularity at the crack tip, with $a_i^J$ and $b_i^{K \alpha}$ being their associated degrees of freedom. Hence, to represent a crack, two sets of additional functions are employed:
\begin{itemize}
\item Heaviside jump function to capture the discontinuity in the displacement. The jump enrichment function is defined as:
\begin{equation}\label{eqn:discontintegrands}
H(\phi) = \left\{ \begin{array}{ccc} 1 & \textup{for} & \phi (x_i) > 0 \\
-1 & \textup{for} & \phi (x_i) < 0
\end{array} \right.
\end{equation}
where $\phi(x_i)$ is the signed distance function from the crack surface defined as:
\begin{equation}
\phi(x_i) = \min_{\overline{x}_i \in \Gamma_c} || x_i - \overline{x}_i|| \textup{sign}( n_i \cdot (x_i - \overline{x}_i ))
\end{equation}
with $n_i$ being the unit outward normal and $\textup{sign}()$ the sign function.
\item Functions with singular derivative near the crack that spans the near tip stress field. For example, in the case of linear elastic fracture mechanics, the following asymptotic displacement field is used:
\begin{equation}\label{eqn:singularintegrands}
F_{\alpha} (r, \theta ) = r^{1/2} \left\{ \sin \frac{\theta}{2}, \cos \frac{\theta}{2}, \sin \frac{\theta}{2} \sin \theta, \cos \frac{\theta}{2} \sin \theta \right\}
\end{equation}
where $r$ is the distance from the crack tip and $\theta$ represents the angular distribution. The linear elastic solution breaks down in the presence of plasticity, with the known HRR fields describing the nature of the dominant singularity instead. Elguedj \textit{et al.} \cite{E06} enriched the shape function basis with the HRR plastic singularity, achieving accurate estimations of standard fracture parameters. Such approach is further extended in this work to incorporate the role of relevant microstructural features (namely, GNDs) in crack tip fields through MSG plasticity. A novel enrichment basis is therefore proposed, where the power of the stress singularity equals $r^{-2/3}$ (see Section \ref{Material model}). As the angular functions play a negligible role in the overall representation of the asymptotic fields \cite{DB08}, the linear elastic fracture mechanics functions are employed.
\end{itemize}

A direct consequence of the enrichment strategy adopted is the possibility of employing simpler meshes that do not need to conform to the crack geometry. The crack can be represented through level sets \cite{D06} or hybrid explicit implicit representation \cite{M11,FB12}. In the present study, a level set representation is used and the enrichment functions at any point of interest are computed using the finite element approximation of the level set functions.

\begin{figure}[H]
\centering
\includegraphics[scale=1.2]{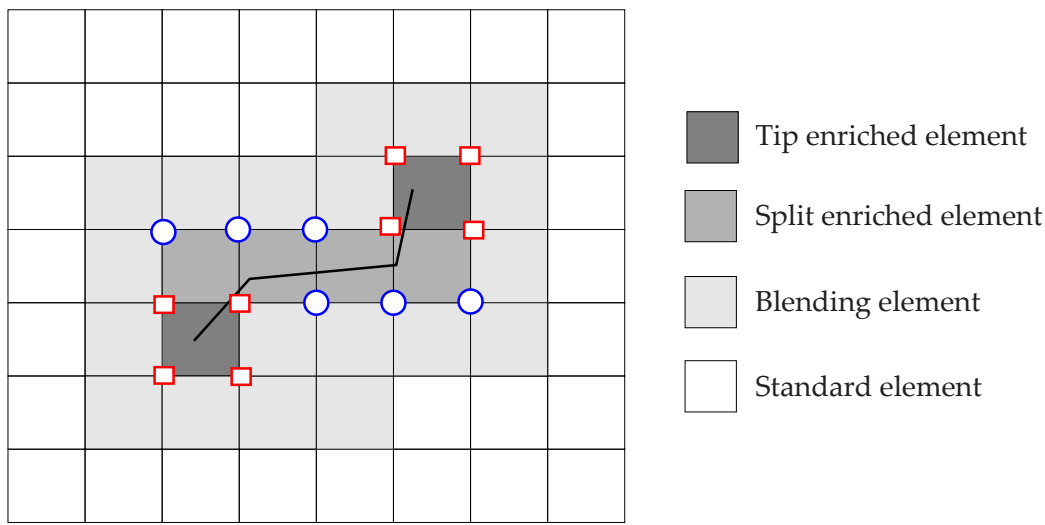}
\caption{Typical X-FEM mesh with an arbitrary crack. Circled nodes are enriched with the discontinuous function while squared nodes are enriched with near-tip asymptotic fields.}
\label{fig:XFEMfig}
\end{figure}

Figure \ref{fig:XFEMfig} shows a typical X-FEM mesh with an arbitrary crack. The enrichment zone is restricted to the vicinity of the crack tip. Elements can be classified into four categories: (a) standard elements; (b) tip enriched elements, (c) split enriched elements and (d) blending elements. The latter are elements at the interface of the standard and enriched elements where the partition of unity is not satisfied and oscillations in the results are observed. This pathological behavior has attracted a considerable research effort and some of the proposed techniques include assumed strain blending elements \cite{RG08}, corrected or weighted XFEM \cite{F08,V08}, hybrid-crack elements \cite{XK07}, semi-analytical approaches \cite{R10,NS13} and spectral functions \cite{L05}. In \cite{LP05}, it was numerically observed that to achieve optimal convergence rate, a fixed area around the crack tip should be enriched with singular functions. This was referred to as \emph{geometrical} enrichment, as opposed to \emph{topological} enrichment, where only one layer of elements around the crack tip is enriched. As detailed below, both topological and geometrical enrichment strategies have been considered in the present work. As proposed by Fries \cite{F08}, a linear weighting function is employed to suppress the oscillatory behavior in the partially enriched elements.\\

Another commonly investigated problem associated with the XFEM is the numerical integration of singular and discontinuous integrands (c.f. Equations \ref{eqn:discontintegrands} - \ref{eqn:singularintegrands}). One potential and yet simple solution for the numerical integration is to partition the elements into triangles. The numerical integration of singular and discontinuous integrands can be alternatively done by: (a) polar integration \cite{LP05}; (b) complex mapping \cite{N09}; (c) equivalent polynomials \cite{V06}; (d) generalized quadrature \cite{MS08}; (e) smoothed XFEM \cite{B11} and (f) adaptive integration schemes \cite{XK06}. Recently, Chin \textit{et al.} \cite{C16} have employed the method of numerical integration of homogeneous functions to integrate discontinuous and weakly singular functions. In the present study, elements are partitioned into triangles and the triangular quadrature rule is employed to integrate the terms in the stiffness matrix.\\

An in-house code is developed in MATLAB for both the FEM and X-FEM cases. Newton-Raphson is employed as solution procedure for the non-linear problem \cite{B09} and stress contours are obtained by performing a Delaunay triangulation and interpolating linearly within the vertex of the triangles (integration points).

\section{Results}
\label{Numerical results}

\subsection{Numerical model}

\begin{figure}[H]
\centering
\includegraphics[scale=0.5]{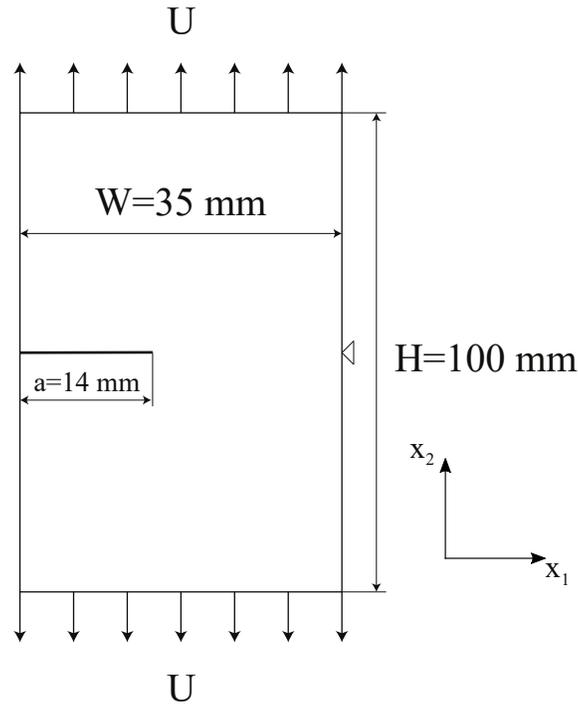}
\caption{Single edge cracked plate: dimensions and boundary conditions.}
\label{fig:Plate}
\end{figure}

As shown in Fig. \ref{fig:Plate}, a cracked plate of dimensions $W=35$ mm (width) and $H=100$ mm (height) subjected to uniaxial displacement is examined. Plane strain conditions are assumed and the horizontal displacement is restricted in the node located at $x_1=W$ and $x_2=H/2$ so as to avoid rigid body motion. The crack is horizontal and located in the middle of the specimen ($H/2$) with the distance from the edge to the tip being $14$ mm. The following material properties are adopted thorough the work: $E=260000$ MPa, $\nu=0.3$, $\sigma_Y=200$ MPa and $n=5$, with isotropic hardening being defined by (\ref{Eq:PowerLaw}). A material length scale of $l=5$ $\mu$m is considered, which would be a typical estimate for nickel \citep{SE98} and corresponds to an intermediate value within the range of experimentally observed material length scales reported in the literature.\\

\subsection{FEM results}

\begin{figure}[H]
\centering
\includegraphics{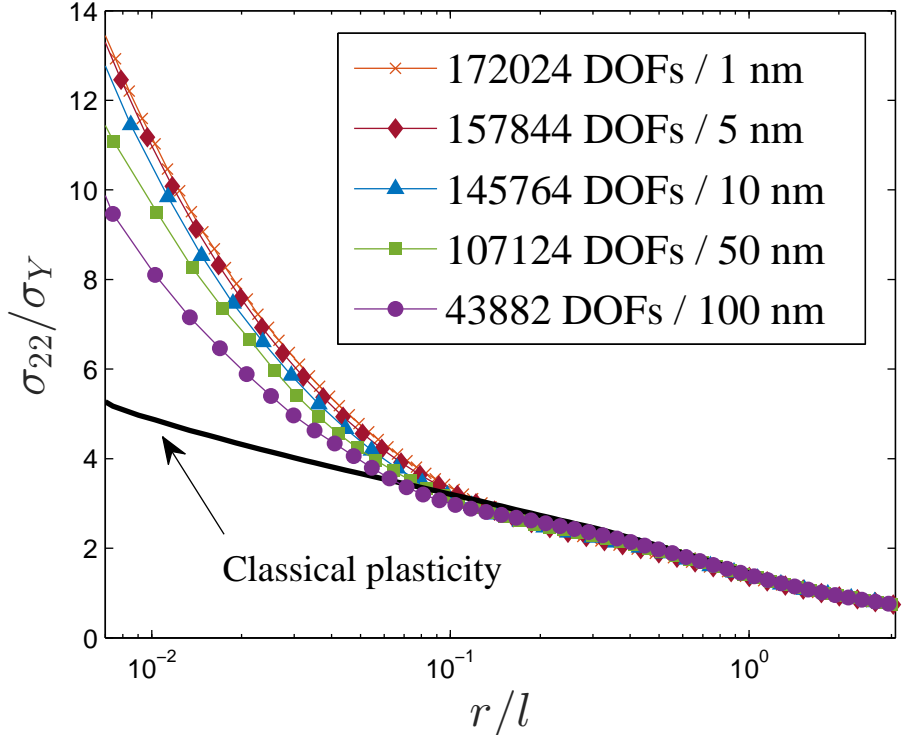}
\caption{Normalized opening stress distribution ahead of the crack tip for different mesh densities, identified as a function of the total number of degrees of freedom (DOFs) and the characteristic length of the element at the crack tip. The figure shows results along the extended crack plane with the normalized distance to the crack tip $r/l$ in log scale.}
\label{fig:FEMresults}
\end{figure}

Fig. \ref{fig:FEMresults} shows the results obtained by means of the standard finite element method for different mesh densities. The legend shows the number of degrees of freedom (DOFs) intrinsic to each mesh, along with the characteristic length of the smallest element in the vicinity of the crack. Quadratic elements with reduced integration have been employed in all cases. The opening stress distribution $\sigma_{22}$ ahead of the crack tip is shown normalized by the initial yield stress while the distance to the crack tip is plotted in logarithmic scale and normalized by the length scale parameter. Results have been obtained for an applied displacement of $U=0.0011$ mm. The prediction obtained for conventional plasticity is also shown in a fine black line and one can easily see that the strain gradient dominated zone is in all cases within $r/l<0.1$ (i.e., 0.5 $\mu$m) for the particular problem, material properties and loading conditions considered.

\begin{figure}[H]
\centering
\includegraphics[scale=0.45]{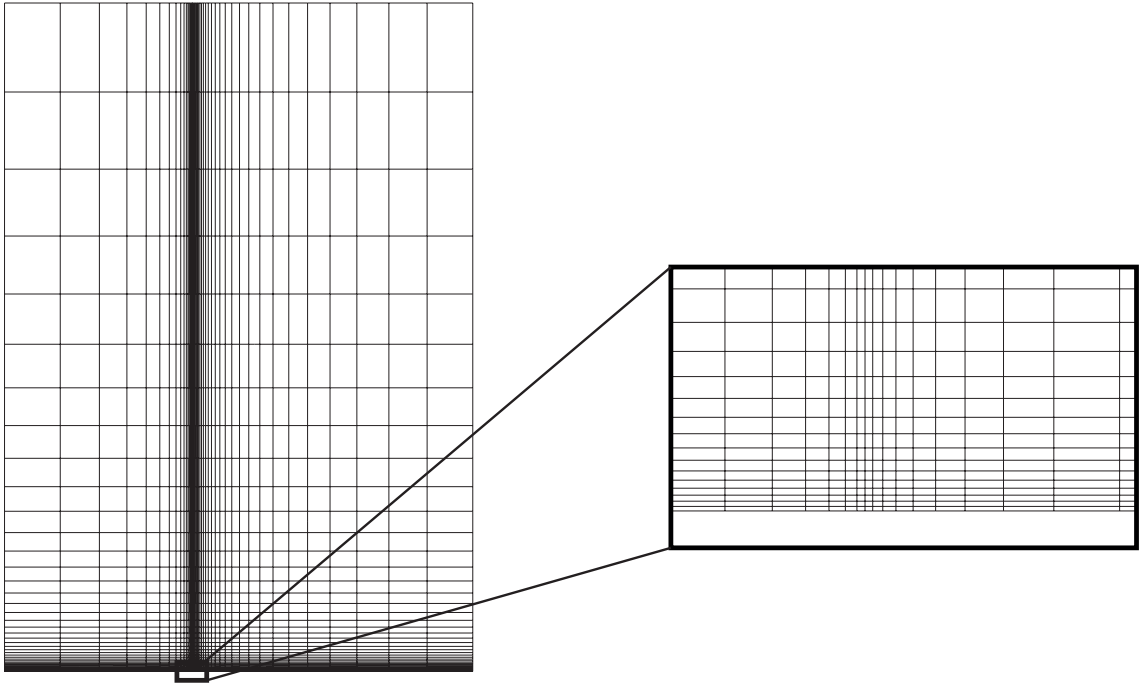}
\caption{Representative finite element mesh, only the upper half of the model is shown due to symmetry.}
\label{fig:FEMmesh}
\end{figure}

Fig. \ref{fig:FEMresults} reveals that numerical convergence has been achieved for a mesh with 157844 DOFs and a characteristic length of the smallest element of 5 nm, as further refinement in the crack tip region leads to almost identical results. This will be considered as the reference finite element solution. A representative illustration of the mesh employed is shown in Fig. \ref{fig:FEMmesh}, where only half of the model is shown, taking advantage of symmetry. As it can be seen in the figure, special care is taken so as to keep an element ratio of 1 close to the crack tip while the mesh gets gradually coarser as we move away from the crack. The use of such small elements is not only very computationally expensive but it also leads to convergence problems as the elements at the crack tip get distorted. Avoiding such level of mesh refinement could strongly benefit fracture and damage assessment within strain gradient plasticity.

\subsection{X-FEM results}

\begin{figure}[H]
\centering
\includegraphics[scale=1.2]{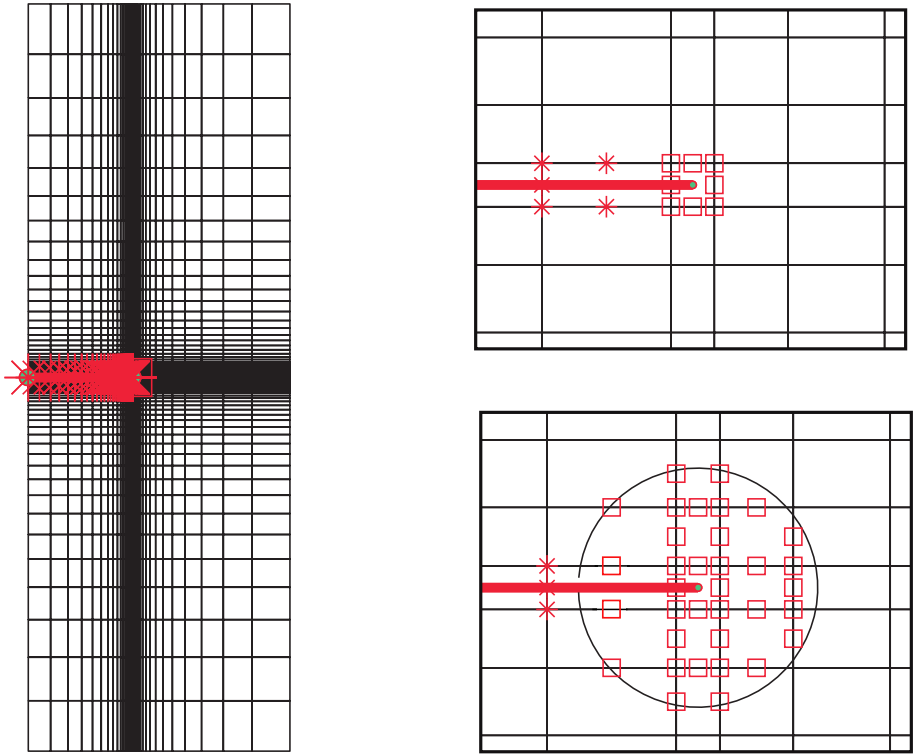}
\caption{Mesh employed in the X-FEM calculations, schematic view and detail of the topological (top) and geometrical (bottom) enrichement regions.}
\label{fig:XFEMmesh1}
\end{figure}

The opening stress is computed in the cracked plate by means of the X-FEM framework described in Section \ref{X-FEM}. A much coarser mesh, relative to the conventional FE case, but with a similar uniform structure is employed, as depicted in Fig. \ref{fig:XFEMmesh1}. A tip element with a characteristic length of 1 $\mu$m is adopted to ensure that the enriched region engulfs the gradient dominated zone. While in the geometrical enrichment case, the characteristic length of the enriched region is chosen so as to coincide with the size of the GNDs-governed domain ($r_e=0.5$ $\mu$m), as discussed below.

\begin{figure}[H]
\centering
\includegraphics[scale=1]{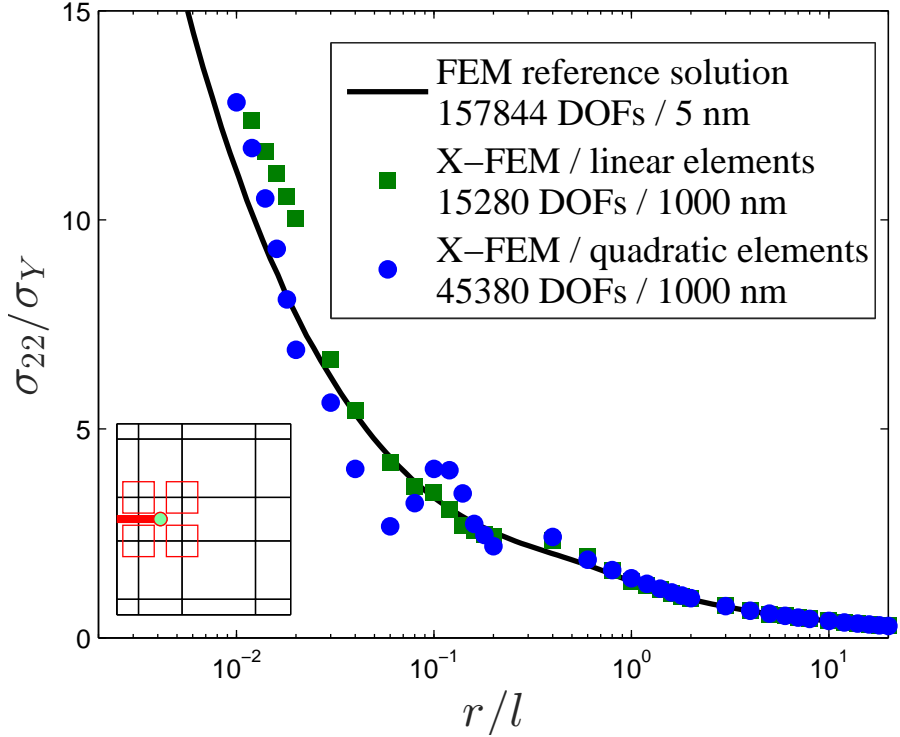}
\caption{Normalized opening stress distribution ahead of the crack tip for topological enrichment, with both linear and quadratic elements, and the reference FEM solution, with  mesh densities identified as a function of the total number of degrees of freedom (DOFs) and the characteristic length of the element at the crack tip. The figure shows results along the extended crack plane with the normalized distance to the crack tip $r/l$ in log scale.}
\label{fig:XFEMresults}
\end{figure}

Results obtained for both quadratic and linear elements are shown in Fig. \ref{fig:XFEMresults}. As in the conventional FE case, the normalized opening stress $\sigma_{22} / \sigma_Y$ is plotted as a function of the normalized distance $r/l$, the latter being in logarithmic scale.\\

X-FEM predictions reveal a good agreement with the reference FE solution, despite the substantial differences in the number of degrees of freedom. Moreover, and unlike the conventional FE case, the influence of strain gradients can also be captured by means of linear quadrilateral elements. This enrichment-enabled capability allows the use of lower order displacement elements, minimizing computational efforts and maximizing user versatility. Further results have been consequently computed with linear elements.

The present gradient-enhanced X-FEM scheme thus shows very good accuracy for a characteristic element length that is two orders of magnitude larger than its standard FEM counterpart. At the local level, small differences are observed in the blending elements, despite the corrected X-FEM approximation adopted. Thereby, enriching the numerical framework with the asymptotic solution of MSG plasticity enables a precise and efficient characterization of crack tip fields, with results being indeed very sensitive to the choice of the power of the stress singularity (see Supplementary Figure 1). Figure \ref{fig:XFEMresultsG} shows the results obtained for a fixed geometrical enrichment radius and different mesh densities.

\begin{figure}[H]
\centering
\includegraphics[scale=1]{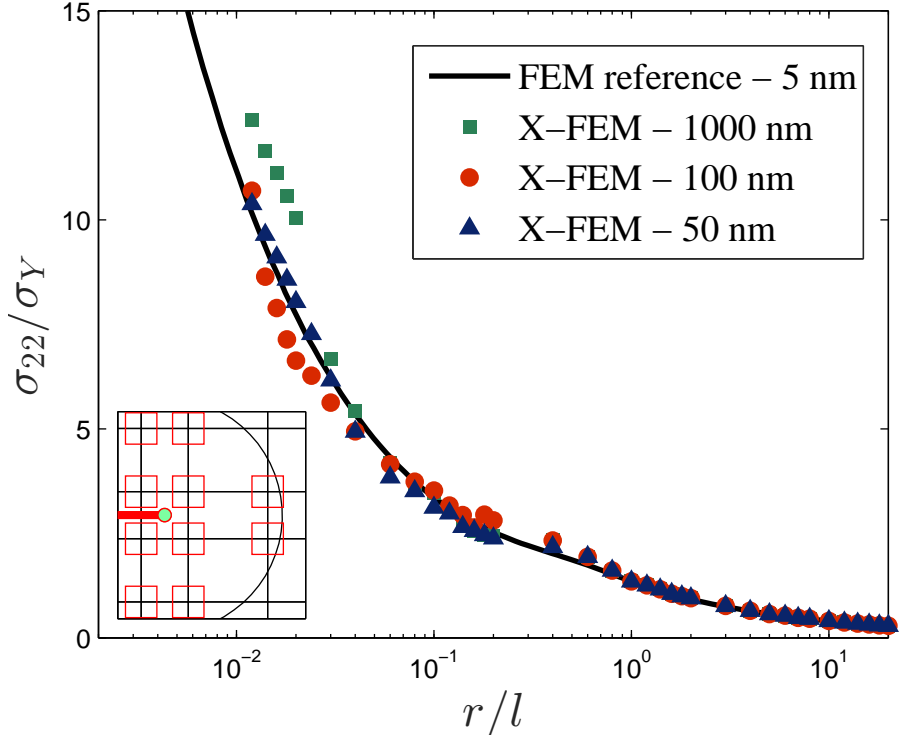}
\caption{Normalized opening stress distribution ahead of the crack tip for geometrical enrichment with enrichment radius $r_e=0.5$ $\mu$m, linear elements and different mesh densities, identified as a function of the characteristic length of the element at the crack tip. The figure shows results along the extended crack plane with the normalized distance to the crack tip $r/l$ in log scale.}
\label{fig:XFEMresultsG}
\end{figure}

As in the topological case, a very promising agreement can be observed, with mesh densities being significantly smaller than the reference FEM solution and computation times varying accordingly. A fixed enrichment radius of $r_e=0.5$ $\mu$m is considered in all cases as the highest precision is achieved when the enriched area and the gradient dominated zone agree. Unlike pure linear elastic analyses, accounting for plastic deformations and the influence of GNDs implies having a crack tip region characterized by three different singular solutions (see Fig. \ref{fig:Domain}). Ideally three classes of asymptotically-enriched nodes should be defined, but such an elaborated scheme is out of the scope of the present work. Hence, the size of the enriched domain must be selected with care to achieve convergence with coarser meshes. This limitation is also intrinsic to the seminal work by Elguedj \textit{et al.} \cite{E06}, where plasticity was confined to the HRR-enriched tip element. The size of the GND-dominated region is nevertheless much less sensitive to material properties or the external load than the plastic zone, and can be properly chosen based on previous parametric studies \cite{MB15}.\\

\begin{figure}[H]
\raggedright
\includegraphics[scale=1.1]{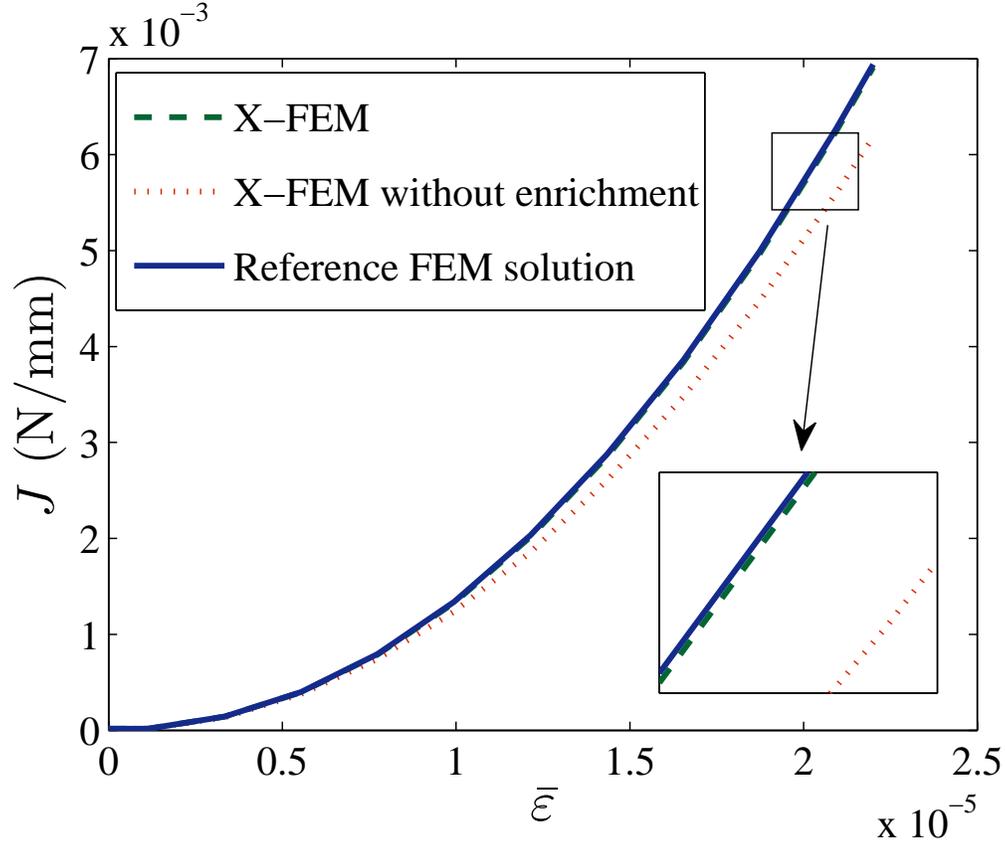}
\caption{$J$-integral versus remote strain for the reference FEM solution ($\ell_e=5$ nm) and the X-FEM solution ($\ell_e=1000$ nm) with and without enrichment. The X-FEM results have been obtained with topological enrichment and linear elements.}
\label{fig:Jresults}
\end{figure}

The capabilities of the proposed numerical scheme to efficiently compute relevant fracture parameters are also examined. First, the $J$-integral is computed by means of the domain integral method \cite{G13} for different load levels, with the external load being characterized through the remote applied strain $\bar{\varepsilon} \equiv 2U/H$. Results obtained are shown in Fig. \ref{fig:Jresults}, where it can be clearly observed that the agreement with the reference FEM solution further increases when a global variable is analyzed. The figure includes the predictions of the standard FEM with a very refined mesh (157844 DOFs and a characteristic element size of $\ell_e=5$ nm) and the results obtained for the present X-FEM scheme, with and without tip enrichment, from a very coarse mesh (15280 DOFs and $\ell_e=1000$ nm). 

\begin{figure}[H]
\centering
\includegraphics[scale=0.9]{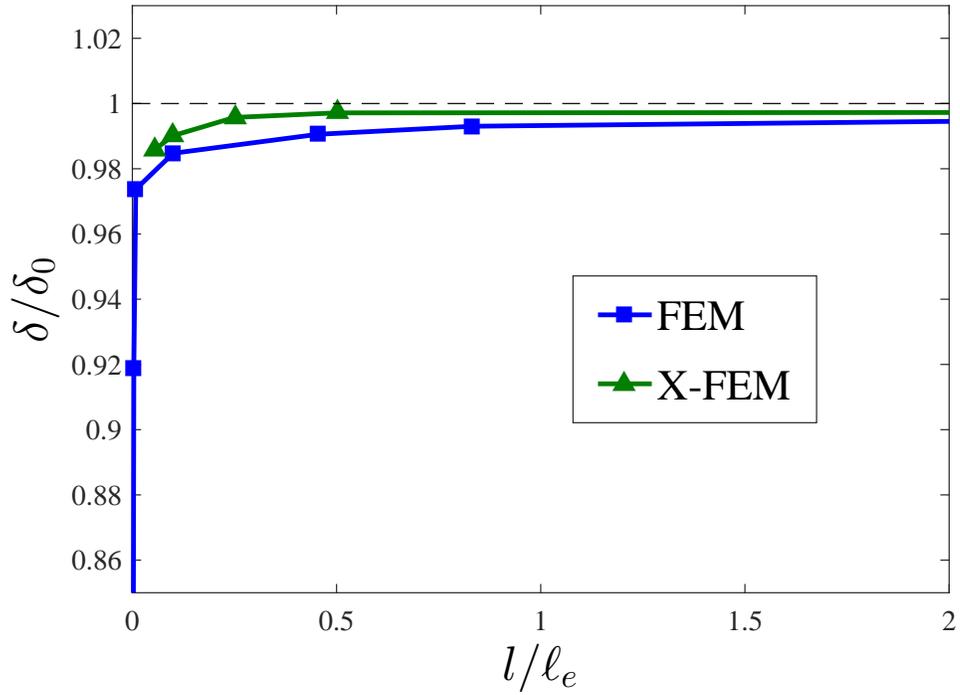}
\caption{Crack opening displacement versus characteristic element size. The horizontal axis corresponds to the material length parameter $l$ divided by the element length $\ell_e$ while the vertical axis is normalized by the crack opening displacement of the reference FEM solution $\delta_0$. The X-FEM results have been obtained with topological enrichment and linear elements.}
\label{fig:ResultsCOD}
\end{figure}

\begin{figure}[H]
\centering
\includegraphics[scale=1]{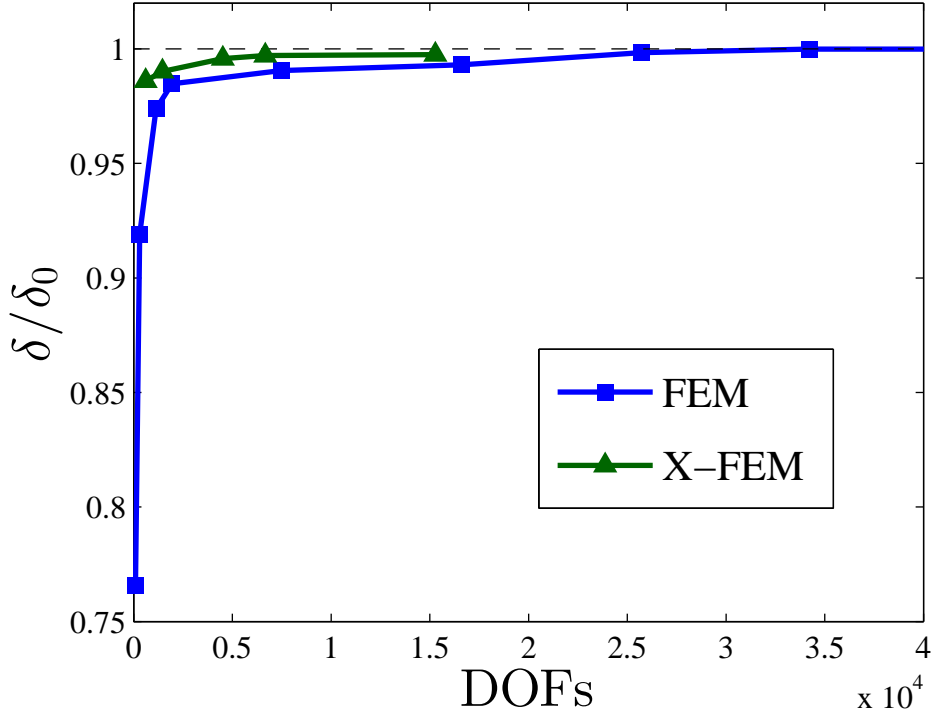}
\caption{Crack opening displacement versus mesh density. The horizontal axis corresponds to the number of degrees of freedom (DOFs) while the vertical axis is normalized by the crack opening displacement of the reference FEM solution $\delta_0$. The X-FEM results have been obtained with topological enrichment and linear elements.}
\label{fig:ResultsCODb}
\end{figure}

The crack opening displacement $\delta$, another meaningful parameter from the fracture mechanics perspective, is also computed for different mesh densities. The magnitude of the crack opening displacement is measured at the crack mouth and its variation with respect to the characteristic element length and the number of degrees of freedom is respectively shown in Figures \ref{fig:ResultsCOD} and \ref{fig:ResultsCODb}. In the former the tip element length is normalized by the length parameter and in both cases the crack opening is shown relative to the reference FEM value.
 
\begin{figure}[H]
\centering
\includegraphics[scale=0.8]{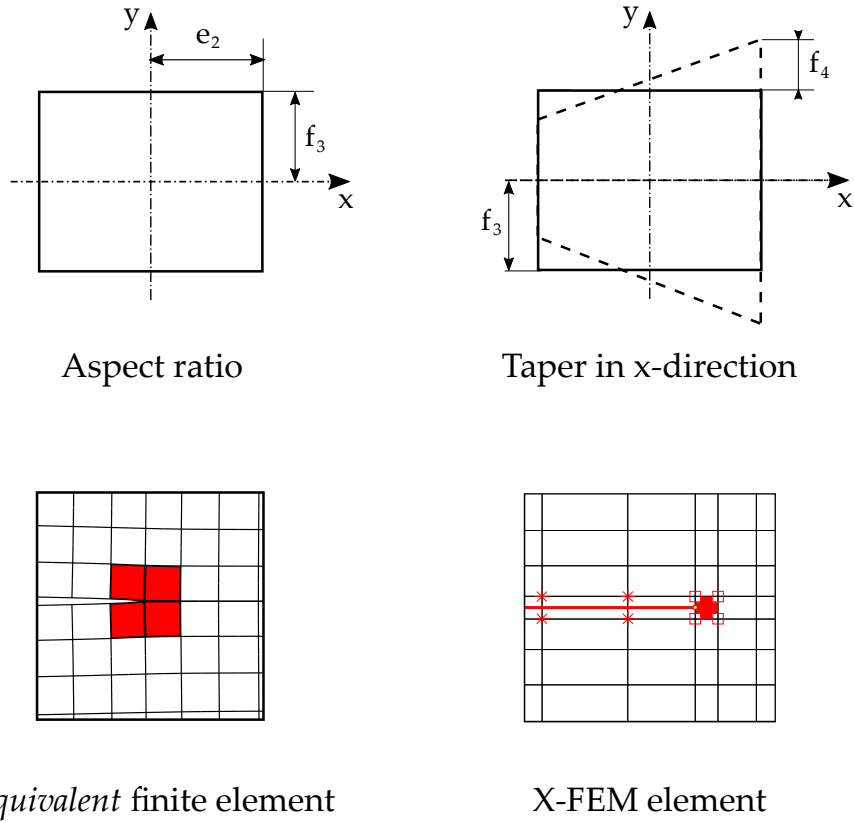}
\caption{Element distortion study: distortion measures (top) and elements under consideration (bottom); the element examined is highlighted in red.}
\label{fig:ElDistortionFig}
\end{figure}

Results reveal a very good performance of the proposed gradient-enhanced enrichment scheme, with an excellent agreement being attained with very coarse meshes. The X-FEM model is able to efficiently track crack tip blunting even when the enriched domain goes far beyond the gradient dominated zone, as the higher stress levels in the conventional plasticity region compensate with the lack of integration points in the vicinity of the crack. The strain gradient plasticity-based enrichment strategy consequently enables accurate estimations of relevant fracture parameters with much less computation effort. This could be of substantial pertinence for structural integrity assessment in engineering industry, where pressures of time and cost demand rapid analyses. As incorporating relevant microstructural features at the micro and nano scales demands intense computations, enriching the FE approximation with the local nature of the solution could be the key enabler for the use of strain gradient plasticity or other multiscale frameworks in crack tip mechanics. Moreover, the present scheme can significantly alleviate the convergence problems intrinsically associated with the use of very refined meshes at deformation levels relevant for practical applications. This has been quantified by investigating the distortion of the crack tip element for both the reference FEM solution and the proposed gradient-enhanced enrichment strategy. As depicted in Fig. \ref{fig:ElDistortionFig} two relevant distortion measures for quadrilateral elements have been considered: the element aspect ratio and the taper in the $x$-direction \cite{R87}; the element skew and the taper in the $y$-direction are zero due to the symmetric nature of the deformation field.

As shown in Fig. \ref{fig:ElDistortionFig} an \emph{equivalent} tip element is defined in the standard FEM approach, so as to directly compare with the X-FEM solution. For this purpose, the crack tip is placed very close to the edge of element in the enriched scheme. By considering only the deformed coordinates of the corner nodes, the following shape parameters are defined:
\begin{equation}
e_2=\frac{1}{4} \left( - x_1 + x_2 + x_3 - x_4 \right)
\end{equation}
\begin{equation}
f_3=\frac{1}{4} \left( - y_1 - y_2 + y_3 + y_4 \right)
\end{equation}
\begin{equation}
f_4=\frac{1}{4} \left( y_1 - y_2 + y_3 - y_4 \right) 
\end{equation}

\noindent where $x_i$ and $y_i$ are respectively the horizontal and vertical \emph{local} nodal coordinates, with counterclockwise node numbering and being the first node the one located in the bottom left corner. The aspect ratio $\Upsilon$ and the taper in the $x$-direction $T^x$ are then defined as:
\begin{equation}
\Upsilon = \textnormal{max} \left\{ \frac{e_2}{f_3} , \frac{f_3}{e_2} \right\}
\end{equation}
\begin{equation}
T^x=f_4/f_3
\end{equation}

\begin{figure}[H]
\centering
\includegraphics[scale=0.8]{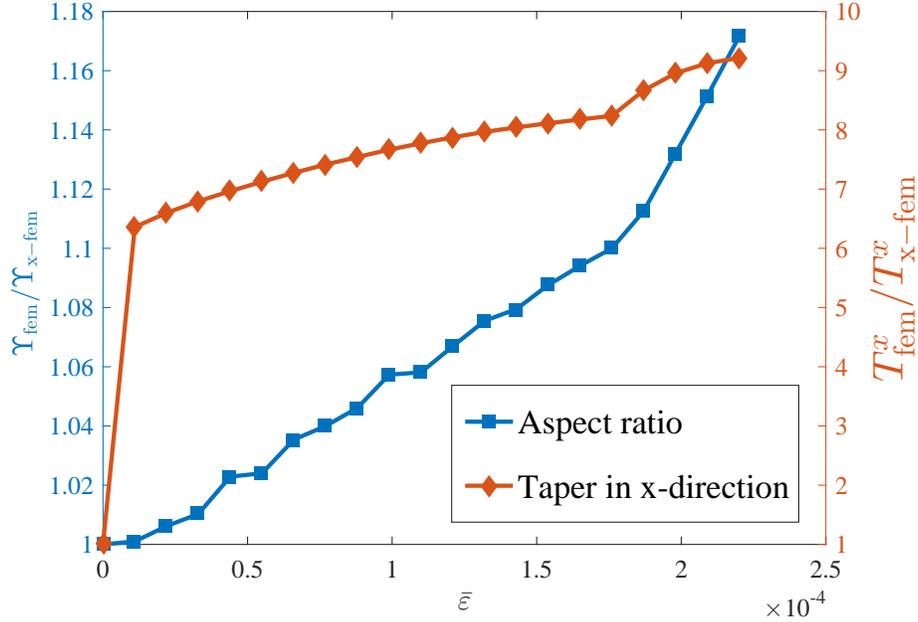}
\caption{Aspect ratio and the taper in the $x$-direction relations between the reference FEM solution and the X-FEM scheme proposed. Results are shown as a function of the external load, characterized by the remote strain $\bar{\varepsilon}$.}
\label{fig:ResultsDist}
\end{figure}

Results obtained as a function of the applied strain are shown in Fig. \ref{fig:ResultsDist}. The predictions for both $\Upsilon$ and $T^x$ with the reference FEM model are shown normalized by the X-FEM results. In that way it can be clearly seen that mesh distortion is severely reduced with an appropriate enriched scheme, even for the relatively low levels of the applied strain considered.  

Significant changes in the element aspect ratio and, particularly, the element tapering take place in the standard FEM model as a consequence of the degree of mesh refinement required. This leads to convergence problems in crack tip mechanics analyses, where elements close to the crack distort excessively \cite{H03,PY11}. The present numerical framework allows to overcome such numerical difficulties and can therefore enable crack tip characterization in a wide range of load levels. This could be particularly useful in environmentally assisted cracking, where GNDs have proven to play a fundamental role \cite{M16a}.

\section{Conclusions}
\label{Concluding remarks}

A robust and efficient numerical framework for crack tip characterization incorporating the role of geometrically necessary dislocations has been developed. The proposed numerical scheme is built from the mechanism-based theory of strain gradient plasticity and takes advantage of its known asymptotic stress singularity to enrich the numerical solution by means of the Extended Finite Element Method (X-FEM). The enriched numerical framework developed can be downloaded from www.empaneda.com/codes and is expected to be helpful to both academic researchers and industry practitioners. The strengths of the proposed X-FEM scheme are clearly seen in the efficient and accurate computation of local stress fields and global fracture parameters; significantly outperforming the standard FEM and avoiding the convergence problems inherent to large element distortions.\\

The range of applicability of the proposed numerical scheme is enormous, as SGP theories have proven to play a fundamental role in a number of structural integrity problems. The use of finite element solutions in large scale engineering applications is hindered by the need to highly refine the mesh in the vicinity of the crack, with the characteristic element length being on the order of a few nanometers. As shown in the present work, this can be readily overcome by employing the X-FEM enriched with the gradient asymptotic solution, SGP-based crack tip characterization being a field where the use of the X-FEM could be of significant relevance. Also, the present numerical framework can be readily develop to model crack propagation within strain gradient plasticity, although a physically-based criterion has yet to be proposed.

\section{Acknowledgments}
\label{Acknowledge of funding}

The authors gratefully acknowledge financial support from the European Research Council (ERC Starting Grant agreement No. 279578). E. Mart\'{\i}nez-Pa\~neda also acknowledges financial support from the Ministry of Economy and Competitiveness of Spain through grant MAT2014-58738-C3-1, and the University of Oviedo through grant UNOV-13-PF.



\end{document}